\documentclass[11pt,leqno,fleqn]{article}
\usepackage{epsf,amsmath,amssymb,amscd,hyperref,srcltx}
\usepackage{amsfonts}
\usepackage{psfrag}
\newcommand{\dps}{\displaystyle}
\newcommand{\End}{\text{\rm End}}
\newcommand{\id}{\text{\rm id}}
\newcommand{\GL}{\text{\rm GL}}
\newcommand{\dy}[2]{%
\refstepcounter{equation}%
\LABEL{#1}%
\begin{list}{}{
\topsep 3mm
\leftmargin 18mm
\rightmargin 0cm
\itemsep 0mm
\listparindent 0mm
\parsep 0mm
\itemsep 0mm
\labelsep 0mm
\labelwidth 18mm
}%
\item[\rm (\theequation)\hfill]
#2
\end{list}%
}
\newcommand{\dyz}[1]{%
\refstepcounter{equation}%
\begin{list}{}{
\topsep 3mm
\leftmargin 18mm
\rightmargin 0cm
\itemsep 0mm
\listparindent 0mm
\parsep 0mm
\itemsep 0mm
\labelsep 0mm
\labelwidth 18mm
}%
\item[\rm (\theequation)\hfill]
#1
\end{list}%
}
\newcommand{\dyyz}[1]{\dyz{\raggedright$\dps#1$}}
\newcommand{\dyy}[2]{\dy{#1}{\raggedright$\dps#2$}}
\newcommand{\di}[2]{%
\refstepcounter{equation}%
\LABEL{#1}%
\begin{list}{}{
\topsep 5mm
\leftmargin 10mm
\rightmargin 0cm
\itemsep 0mm
\listparindent 0mm
\parsep 0mm
\labelsep 1mm
\labelwidth 10mm
}%
\item[\rm (\theequation)\hfill]
\begin{list}{}{
\topsep 0mm
\leftmargin 8mm
\rightmargin 0mm
\itemsep 0mm
\listparindent 0mm
\parsep 0mm
\labelsep 1.5mm
\labelwidth 6.5mm
}
#2
\end{list}%
\end{list}%
}
\newcommand{\nr}[1]{\item[{\rm (#1)}]}
\newcommand{\nrs}[1]{\item[{\rm (#1)}]\vspace{-\itemsep}}
\newcounter{stelling}

\newcommand{\thmnn}[1]{\setcounter{claim}{0}\vspace{4mm}\noindent{\bf Theorem.}{\it #1}}
\newcounter{claim}
\newcommand{\cl}[2]{\refstepcounter{claim}\vspace{4mm}\noindent{\em Claim \theclaim.} \label{#1} {\it #2}}
\newcommand{\clz}[1]{\refstepcounter{claim}\vspace{4mm}\noindent{\em Claim \theclaim.}  {\it #1}}
\newcounter{sectie}
\newcommand{\sect}[2]{\refstepcounter{sectie}
\section*{\boldmath \thesectie. #2}%
\label{#1}}
\newcommand{\sectz}[1]{\refstepcounter{sectie}
\section*{\boldmath \thesectie. #1}%
}
\newcommand{\pf}{\vspace{3mm}\noindent{\bf Proof.}\ }
\newcommand{\pfcl}{\vspace{3mm}\noindent{\em Proof.}\ }
\newcommand{\bx}{\hspace*{\fill} \hbox{\hskip 1pt \vrule width 4pt height 8pt depth 1.5pt \hskip 1pt}

\addvspace{4mm}}
\newcommand{\obx}{\hspace*{\fill} \hbox{$\Box$}

\addvspace{4mm}}
\newcommand{\openbx}{\hspace*{\fill} $\Box$\\ \vspace{1mm}}
\newcommand{\rf}[1]{{\rm (\ref{#1})}}

\newcommand{\tr}{\text{\rm tr}}
\newcommand{\AAA}{{\cal A}}
\newcommand{\CC}{{\cal C}}

\newcommand{\II}{{\cal I}}
\newcommand{\OO}{{\cal O}}
\newcommand{\TT}{{\cal T}}
\newcommand{\VV}{{\cal V}}
\newcommand{\XX}{{\cal X}}
\newcommand{\LABEL}[1]{\label{#1}}
\newcommand{\rank}{\text{\rm rank}}
\newcommand{\sgn}{\text{\rm sgn}}
\newcommand{\oC}{{\mathbb{C}}}
\newcommand{\oR}{{\mathbb{R}}}
\newcommand{\oZ}{{\mathbb{Z}}}
\renewcommand{\loop}{\hspace*{-6pt}\raisebox{-.2\height}{\scalebox{0.05}{\includegraphics{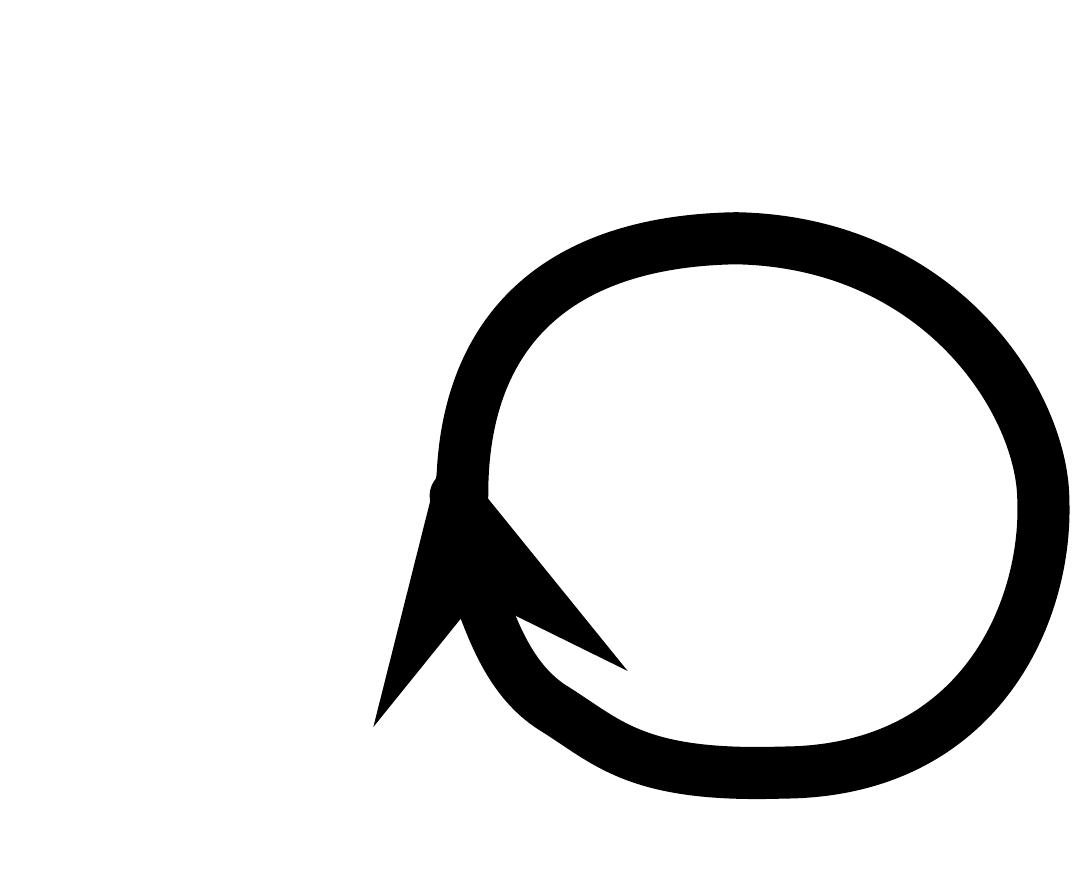}}}}
\newcommand{\Ker}{{\text{\rm Ker}}}
\newcommand{\gl}{{\frak{gl}}}
\newcommand{\tussenkop}[1]{\medskip
\noindent
{\bf #1.}}
\begin{document}

\begin{center}
{\large\bf CONNECTION MATRICES AND LIE ALGEBRA WEIGHT SYSTEMS FOR MULTILOOP CHORD DIAGRAMS

}
\vspace{4mm}

{\large
Alexander Schrijver\footnote{ University of Amsterdam and CWI, Amsterdam.
Mailing address: CWI, Science Park 123, 1098 XG Amsterdam,
The Netherlands.
Email: lex@cwi.nl.
The research leading to these results has received funding from the European Research Council
under the European Union's Seventh Framework Programme (FP7/2007-2013) / ERC grant agreement
n$\mbox{}^{\circ}$ 339109.}}

\end{center}

\noindent
{\small
{\bf Abstract.}
We give necessary and sufficient conditions for a weight system on multiloop chord diagrams to be
obtainable from a metrized Lie algebra representation, in terms of a bound on the ranks of
associated connection matrices.

Here a multiloop chord diagram is a graph with directed and undirected edges so that at each vertex
precisely one directed edge is entering and precisely one directed edge is leaving, and each
vertex is incident with precisely one undirected edge.
Weight systems on multiloop chord diagrams yield the Vassiliev invariants for knots and links.

The $k$-th connection matrix of a function $f$ on the collection of multiloop chord diagrams is
the matrix with rows and columns indexed by $k$-labeled chord tangles, and with entries equal to the
$f$-value on the join of the tangles.

}

\sectz{Introduction}

In this introduction we describe our results for those familiar with the basic theory of
weight systems on chord diagrams (cf.\ [4]).
In the next section we define concepts, so as to fix terminology and so as to 
make the paper self-contained also for those not familiar with weight systems.

Bar-Natan [1,\linebreak[0]2] and Kontsevich [10] have shown that any finite-dimensional
representation $\rho$ of a metrized Lie algebra $\frak{g}$ yields a weight system
$\varphi_{\frak{g}}^{\rho}$ on chord diagrams --- more generally, on multiloop chord diagrams.
(These are chord diagrams in which more than one Wilson loop is allowed.
Weight systems on multiloop chord diagrams yield Vassiliev {\em link} invariants.)

In this paper, we characterize the weight systems that arise this way.
More precisely, we show the equivalence of the following conditions for any
complex-valued weight system $f$:
\dy{14se14d}{
(i) $f=\varphi_{\frak{g}}^{\rho}$ for some completely reducible faithful representation $\rho$ of some metrized Lie algebra $\frak{g}$;\\
(ii) $f=\varphi_{\frak{g}}^{\rho}$ for some representation $\rho$ of some metrized Lie algebra $\frak{g}$;\\
(iii) $f$ is the partition function $p_R$ of some $n\in\oZ_+$ and $R\in S^2(\gl(n))$;\\
(iv) $f(\loop)\in\oR$ and $\rank(M_{f,k})\leq f(\loop)^{2k}$ for each $k$.
}
Throughout, $\gl(n)=\gl(n,\oC)$, while $\oC$ may be replaced by any algebraically closed field of characteristic 0.
All representations are assumed to be finite-dimensional.
In (i), the Lie algebra $\frak{g}$ is necessarily reductive.
The largest part of the proof consists of showing (iv)$\Longrightarrow$(iii).

We give some explanation of the conditions (iii) and (iv).
First, $S^2(\gl(n))$ denotes the space of tensors in $\gl(n)\otimes\gl(n)$ that are symmetric
(i.e., invariant under the linear function induced by $X\otimes Y\to Y\otimes X$).
The {\em partition function} $p_R$ of $R\in S^2(\gl(n))$ can be intuively described as
the function on multiloop
chord diagrams obtained by inserting a copy of the tensor $R$ at each chord,
assigning (`multilinearly') its two tensor components in $\gl(n)$ to the two ends of that chord,
next calculating, along any Wilson loop, the trace of
the product of the elements in $\gl(n)$ assigned to the vertices of that Wilson loop (in order),
and finally taking the product of these traces over all Wilson loops.
(This is in analogy to the partition function of the `vertex model' in de la Harpe and Jones [8].)

In \rf{14se14d}(iv), $\loop$ is the chord diagram without chords.
To describe the matrix $M_{f,k}$, we need `$k$-labeled multiloop chord tangles', or `$k$-tangles' for short.
A {\em $k$-tangle} is a multiloop chord diagram with $k$ directed edges entering it, labeled
$1,\ldots,k$, and $k$ directed edges leaving it, also labeled $1,\ldots,k$, like the 4-tangle
\begin{center}
\raisebox{-.2\height}{\scalebox{0.1}{\includegraphics{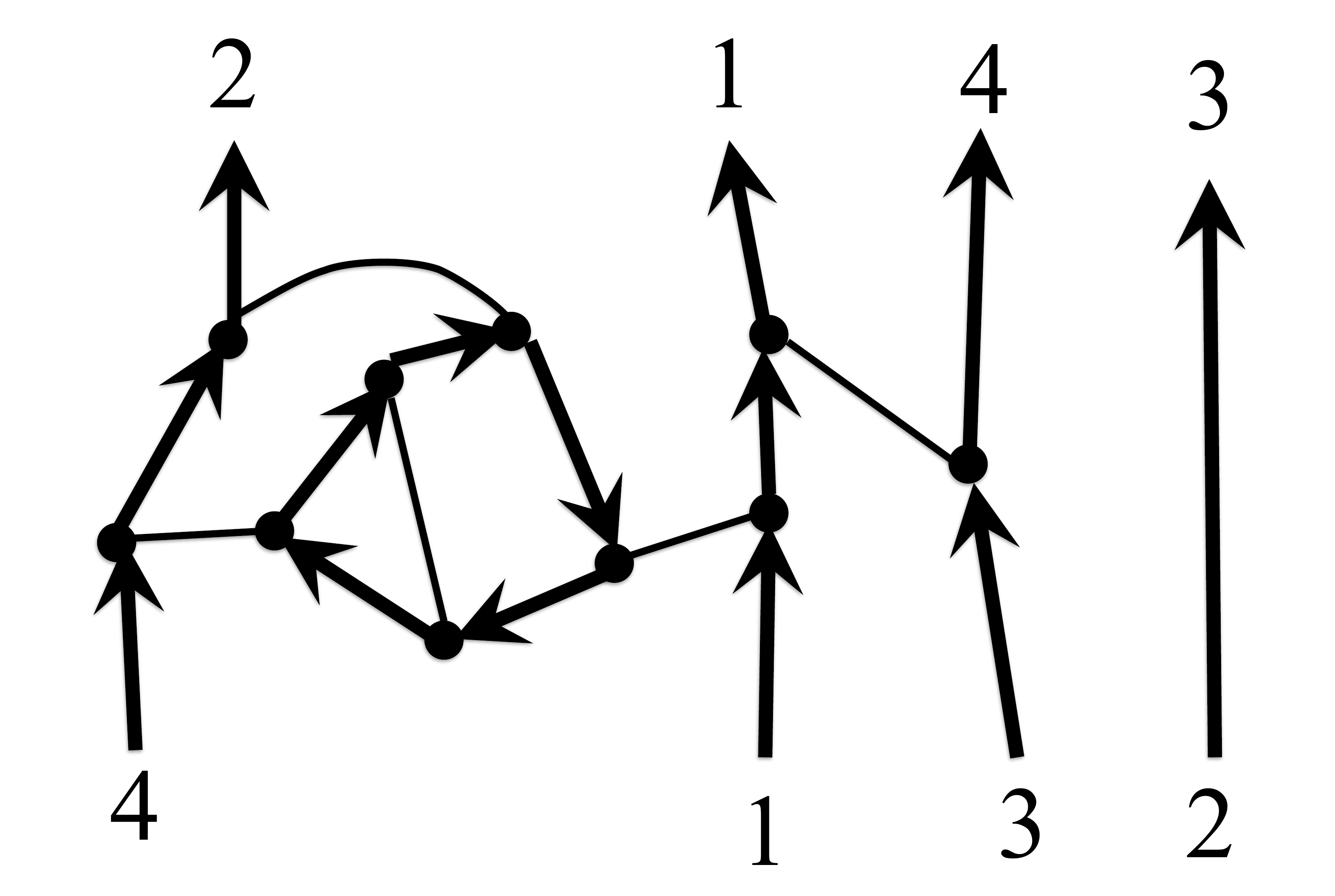}}}.
\end{center}
Let $\TT_k$ denote the collection of all $k$-tangles.
For $S,T\in\TT_k$, let $S\cdot T$ be the multiloop chord diagram
obtained by glueing $S$ and $T$ appropriately together:
$S\cdot T$ arises from the disjoint union of $S$ and $T$ by identifying
outgoing edge labeled $i$ of $S$ with ingoing edge labeled $i$ of $T$, and similarly,
identifying
outgoing edge labeled $i$ of $T$ with ingoing edge labeled $i$ of $S$
(for $i=1,\ldots,k$).
Then the {\em $k$-th connection matrix} $M_{f,k}$ of $f$ is the $\TT_k\times\TT_k$ matrix with entry $f(S\cdot T)$ in position
$(S,T)\in\TT_k\times\TT_k$.
(Studying such matrices roots in work of Freedman, Lov\'asz, and Schrijver [6] and Szegedy [14],
cf.\ also the recent book by Lov\'asz [12].)

The implications (i)$\Longrightarrow$(ii)$\Longrightarrow$(iii)$\Longrightarrow$(iv)
are easy --- the content of this paper is proving the reverse implications.
Indeed, (i)$\Longrightarrow$(ii) is trivial.
To see (ii)$\Longrightarrow$(iii), recall the fundamental construction of Bar-Natan [1,\linebreak[0]2] and
Kontsevich [10].
Let $\frak{g}$ be a metrized Lie algebra and let $\rho:\frak{g}\to\gl(n)$ be a representation.
Let $b_1,\ldots,b_k$ be any orthonormal basis of $\frak{g}$ and define
\dyyz{
R(\frak{g},\rho):=\sum_{i=1}^k\rho(b_i)\otimes\rho(b_i)\in S^2(\gl(n))
}
(which is independent of the choice of the orthonormal basis).
Then $\phi_{\frak{g}}^{\rho}:= p_{R(\frak{g},\rho)}$ is a weight system.
So one has (ii)$\Longrightarrow$(iii).

The Lie bracket is not involved in condition \rf{14se14d}(iii), it is required only that
$p_R$ be a weight system.
Indeed, not each $R\in S^2(\gl(n))$ for which $p_R$ is a weight system arises as above from a Lie algebra.
For instance, let
\arraycolsep 2pt
$B_1:=$
{\tiny
$\left(
\begin{array}{cc}
1&1\\
0&0
\end{array}
\right)$
}
and
$B_2:=$
{\tiny
$\left(
\begin{array}{cc}
0&1\\
0&1
\end{array}
\right)$
}
(as elements of $\gl(2)$),
and set $R:=B_1^{\otimes 2}+B_2^{\otimes 2}\in S^2(\gl(2))$.
Then $p_R$ is identically 2 on connected diagrams,
hence $p_R$ is a weight system, but there is no representation $\rho$ of a
metrized Lie algebra $\frak{g}$ with $R=R(\frak{g},\rho)$
(essentially because the matrices $B_1$ and $B_2$ do not span a matrix Lie algebra).

The implication (iii)$\Longrightarrow$(iv) follows from the fact that for any $k$ and any
$k$-tangles $S$ and $T$, $p_R(S\cdot T)$ can be described as
the trace of the product of certain elements $\widehat p_R(S)$ and $\widehat p_R(T)$ of $\gl(n)^{\otimes k}$,
where the latter space has dimension $n^{2k}$.

Our proof of the reverse implications is based on some basic results of algebraic geometry (Nullstellensatz),
invariant theory (first and second fundamental theorem, closed orbit theorem), and
(implicitly through [13]) the representation theory of the symmetric group.
It consists of showing that if \rf{14se14d}(iv) is satisfied, then $n:=f(\loop)$ belongs to
$\oZ_+$ and the affine $\GL(n)$-variety
\dyyz{
\VV:=\{R\in S^2(\gl(n))\mid p_R=f\}
}
is nonempty (which is (iii)), and each $R$ in the (unique) closed $\GL(n)$-orbit in $\VV$ produces a
completely reducible faithful representation of a Lie algebra as in (i).

We must emphasize here that the above will be proved for {\em multiloop} chord diagrams.
We do not know in how far it remains true when restricting the functions to ordinary, one-loop,
chord diagrams.

We also do not know in how far the Lie algebra $\frak{g}$ and the representation $\rho$ in \rf{14se14d}(i) are unique
(up to the action of $\GL(n)$ where $n$ is the dimension of $\rho$),
although the existence is shown by construction from the unique closed $\GL(n)$-orbit in $\VV$.
A partial result in this direction was given by Kodiyalam and Raghavan [9]:
let $\frak{g}$ and $\frak{g}'$ be $n$-dimensional semisimple Lie algebras,
with the Killing forms as metrics,
and let $\rho$ and $\rho'$ be the adjoint representations;
if $\varphi_{\frak{g}}^{\rho}=\varphi_{\frak{g}'}^{\rho'}$ on (one-loop) chord diagrams,
then $\frak{g}=\frak{g}'$.

\sect{22jn14a}{Preliminaries}

\tussenkop{Multiloop chord diagrams and weight systems}
A {\em multiloop chord diagram} is a cubic graph $C$ in which a collection
of disjoint oriented cycles is specified that cover all vertices.
These cycles are called the {\em Wilson loops}, and the remaining edges
(that form a perfect matching on the vertex set of $C$) are called the {\em chords}.

Alternatively, a multiloop chord diagram can be described as a graph with directed
and undirected edges such that for each vertex $v$:
\dy{22me14a}{
$v$ is entered by precisely one
directed edge, is left by precisely one directed edge, and is incident
with precisely one undirected edge,
}
as in
\raisebox{-.5\height}{\scalebox{0.1}{\includegraphics{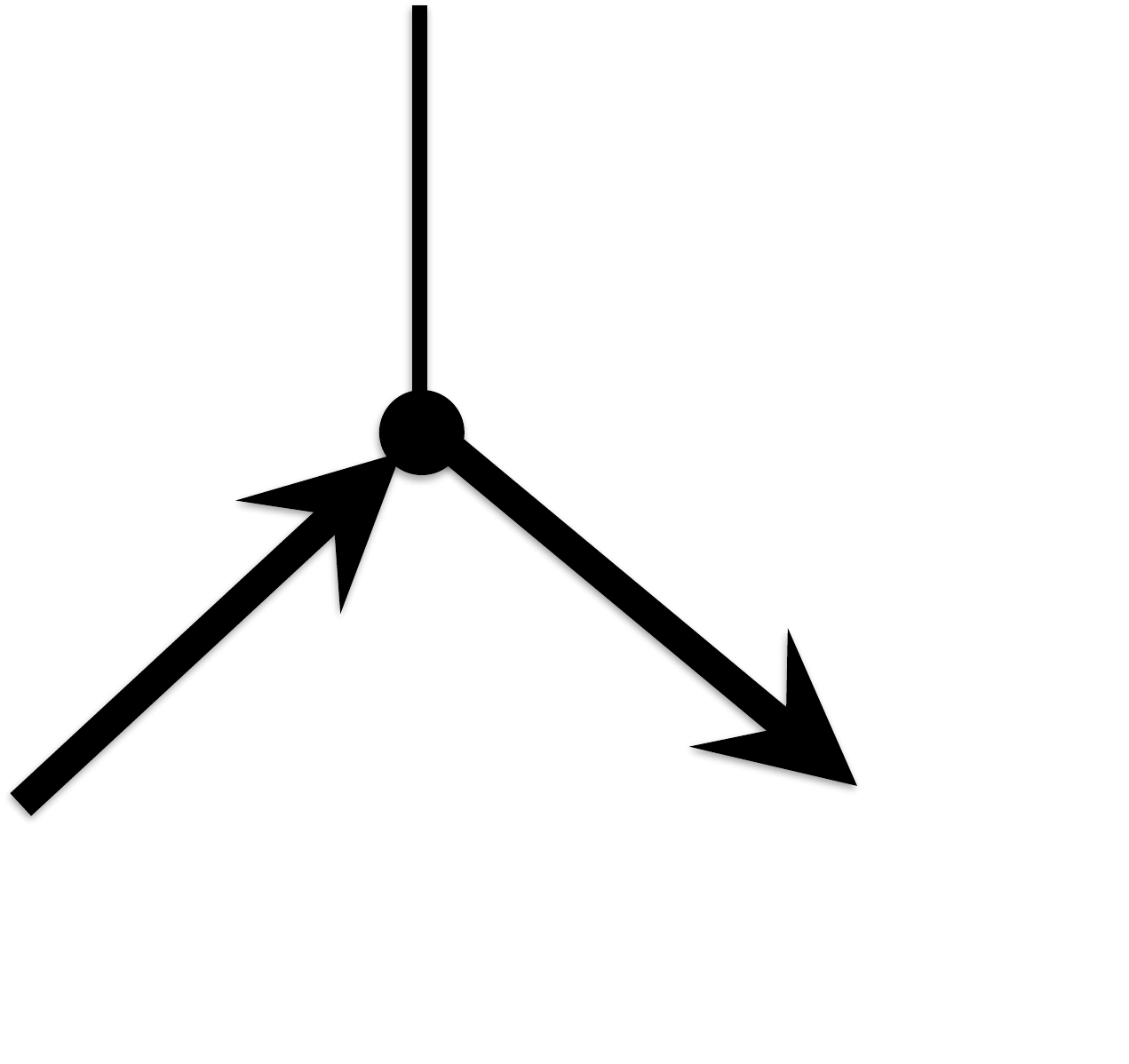}}}.
The following is an example of a multiloop chord diagram:
\vspace*{-8mm}
\begin{center}
\raisebox{-.2\height}{\scalebox{0.3}{\includegraphics{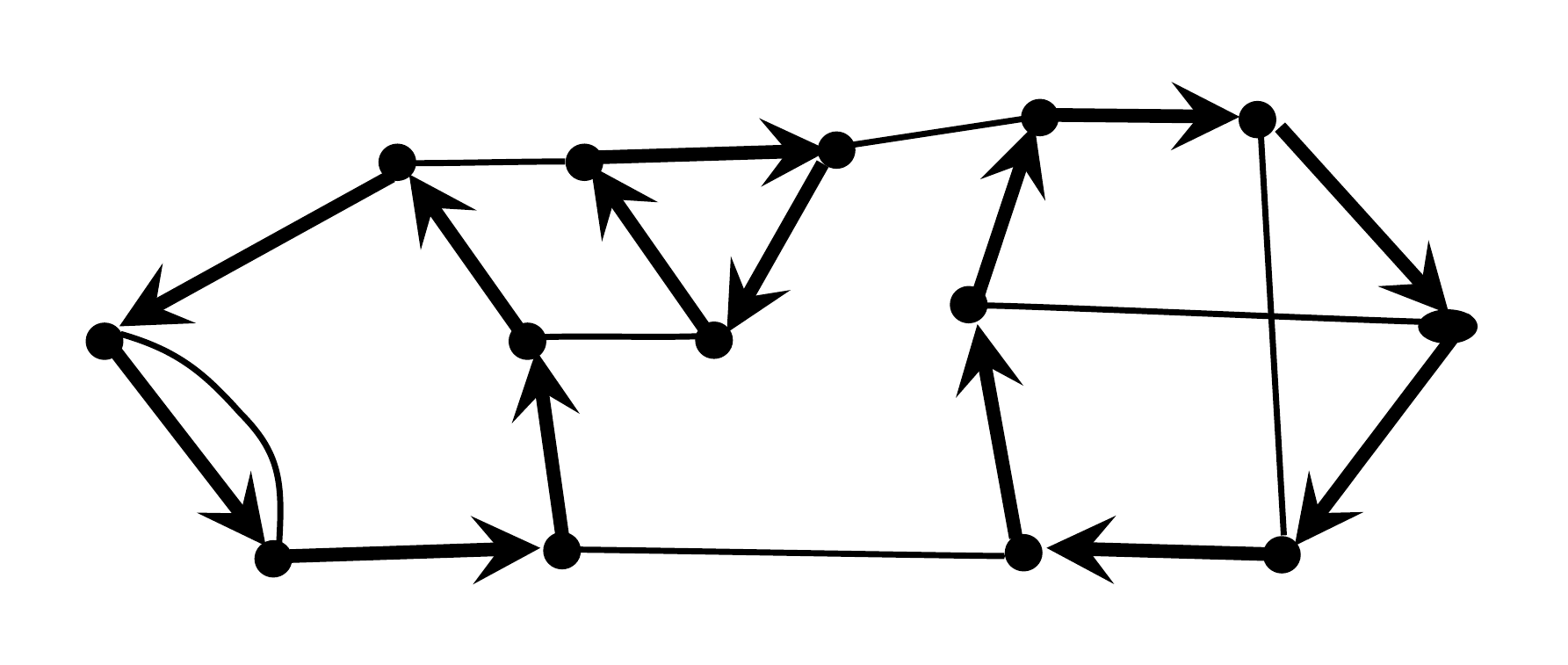}}}.
\end{center}
\vspace*{-3mm}
Directed loops are allowed, but no undirected loops.
Moreover, we allow the `vertexless directed loop' $\loop$ (in other words, the chord diagram
of order 0) --- more precisely,
components of a multiloop chord diagram may be vertexless directed loops.

Let $\CC$ denote the collection of multiloop chord diagrams.
Basic for Vassiliev knot invariants
(cf.\ [4]) are functions $f$ on $\CC$ that satisfy certain linear relations, called
the {\em 4-term (4T) relations}.
They can be visualized as:
\begin{center}
$
f\big(\raisebox{-.2\height}{\scalebox{0.078}{\includegraphics{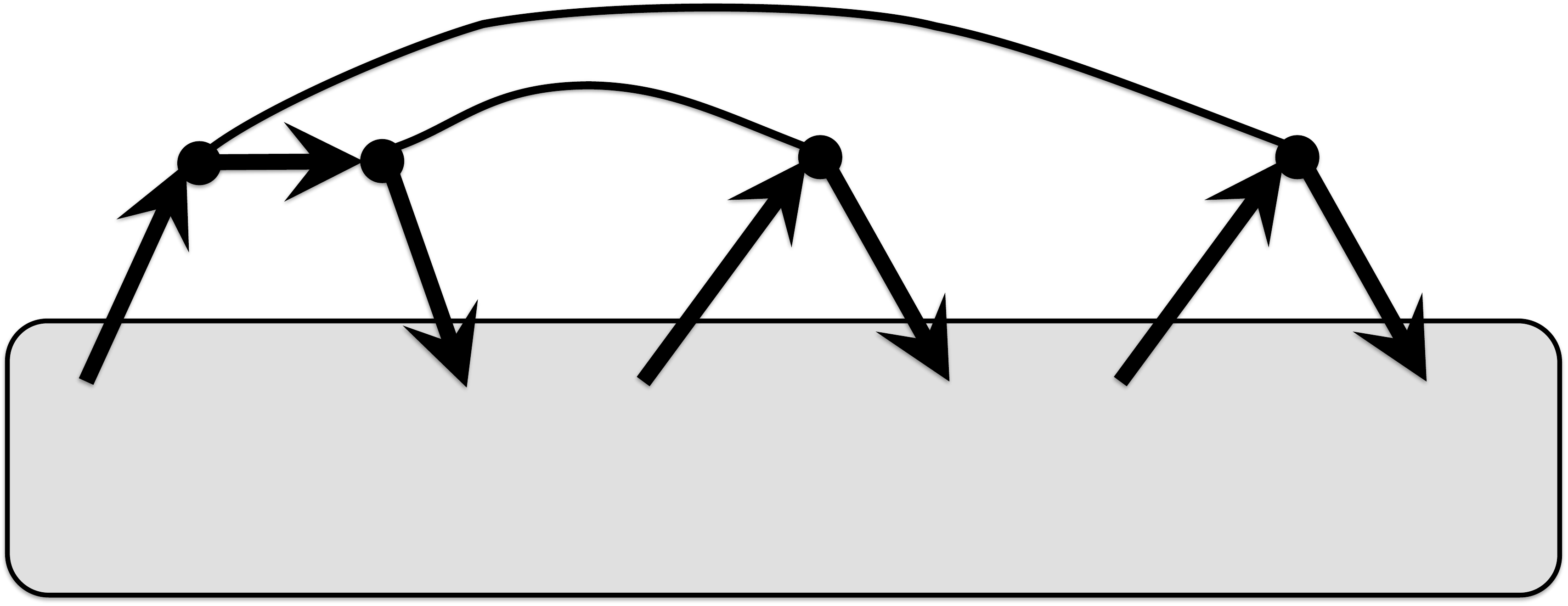}}}\big)
-
f\big(\raisebox{-.2\height}{\scalebox{0.078}{\includegraphics{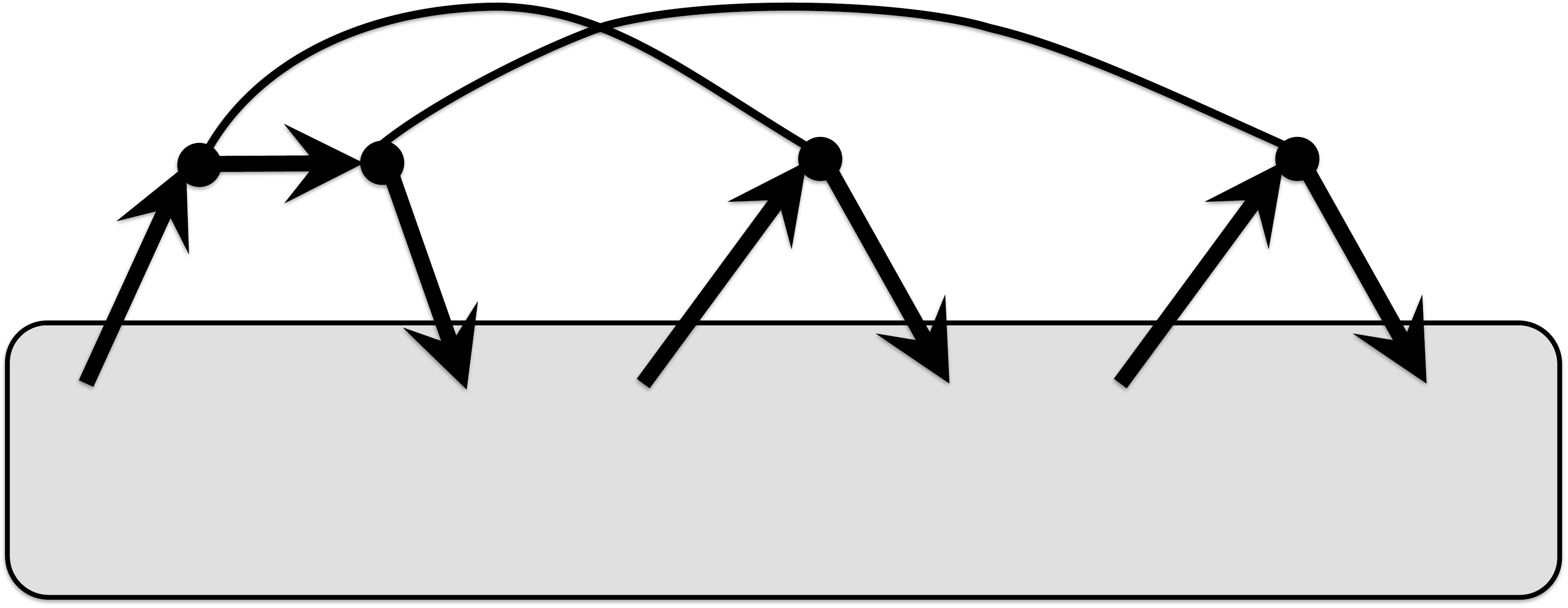}}}\big)
=
f\big(\raisebox{-.2\height}{\scalebox{0.078}{\includegraphics{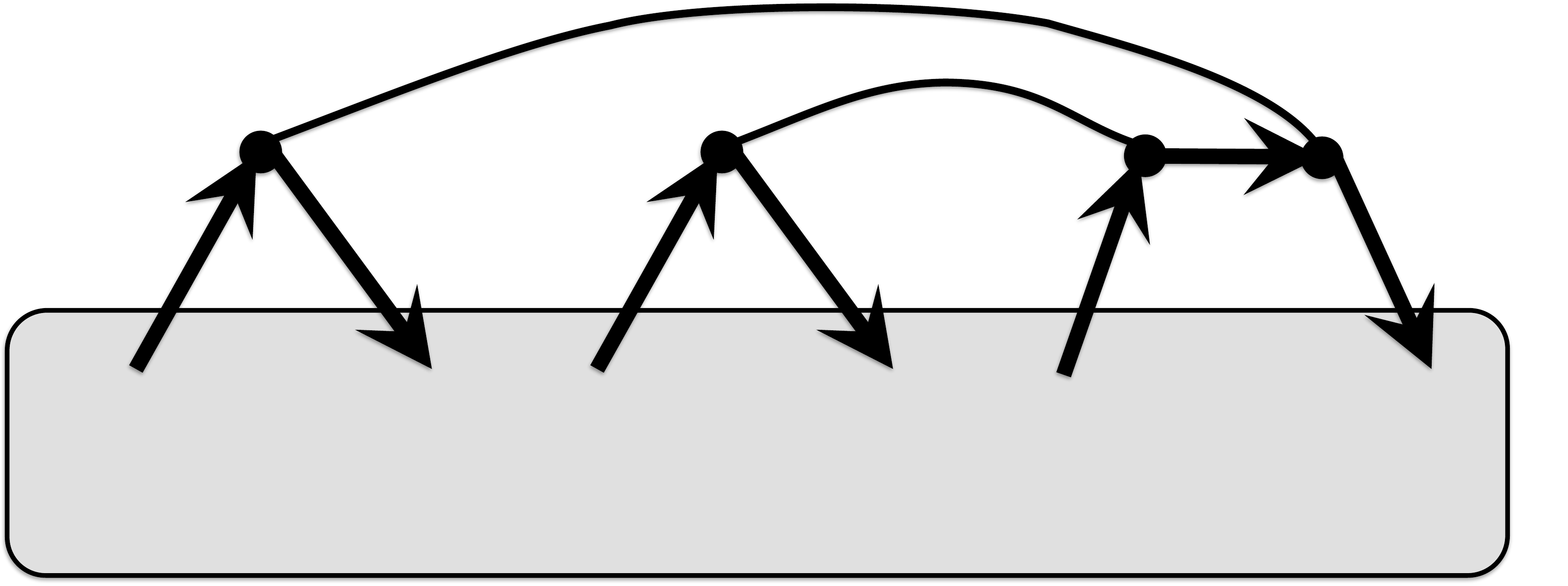}}}\big)
-
f\big(\raisebox{-.2\height}{\scalebox{0.078}{\includegraphics{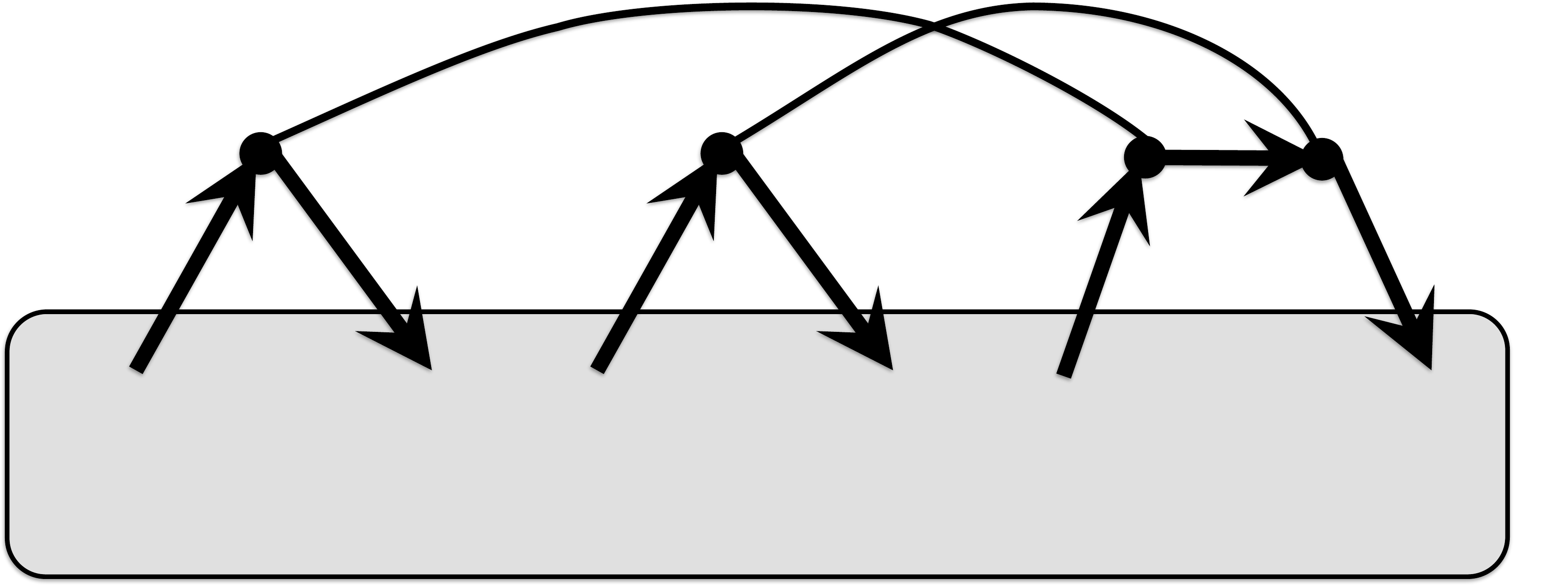}}}\big).
$
\end{center}
(Each of the four grey rectangles contains the rest of the diagram, the same in each rectangle.)

Functions satisfying the 4T relations are called {\em weight systems}.
More precisely, we call a function $f$ on multiloop chord diagrams a {\em weight system} if it satisfies the
4T relations, and moreover it is {\em multiplicative}: $f(\emptyset)=1$ and $f(C\sqcup D)=f(C)f(D)$
for all multiloop chord diagrams $C,D$, where $C\sqcup D$ denotes the disjoint union of
$C$ and $D$.
Hence any weight system is determined by its values on connected multiloop chord diagrams.

Through the {\em Kontsevich integral}, each $\oC$-valued weight system on the collection of multiloop chord diagram
with some fixed number of chords and some fixed number $t$ of Wilson loops, gives an invariant
for links with $t$ components.
They produce precisely the {\em Vassiliev invariants} for knots and links.
We refer for these important concepts to the book of Chmutov, Duzhin, and Mostovoy [4] --- for understanding our treatment below
they are however not needed.

\tussenkop{Some notation and linear algebra}
As usual,
\dyz{
$\oZ_+:=$ the set of nonnegative integers, and $[k]:=\{1,\ldots,k\}$
}
for any $k\in\oZ_+$.

For any set $\XX$, $\oC\XX$ denotes the linear space
of formal $\oC$-linear combinations of finitely many elements of $\XX$.
(Occasionnally, elements of $\oC\XX$ are called {\em quantum} elements of $\XX$.)
Any function on $\XX$ to a $\oC$-linear space can be uniquely extended to a linear function on $\oC\XX$.

For a linear space $X$, $S^2(X)$ denotes the space of symmetric elements of $X\otimes X$, i.e.,
those invariant under the linear operation on $X\otimes X$ induced by $x\otimes y\mapsto y\otimes x$.
It is elementary matrix theory to prove that if $X$ is finite-dimensional,
then for any $R\in S^2(X)$ there is a unique subspace
$Y$ of $X$ and a unique nondegenerate bilinear form on $Y$ such that
for each orthonormal basis $b_1,\ldots,b_k$ of $Y$ one has
\dyy{4se14f}{
R=\sum_{i=1}^kb_i\otimes b_i.
}
Considering $R$ as matrix, $Y$ is equal to the column space of $R$.

\tussenkop{Partition functions}
Each $R\in S^2(\gl(n))$ gives a function $p_R$ on multiloop chord diagrams as follows.
Fix a basis of $\oC^n$, and write $R=(R_{i,j}^{k,l})$, with $i,j,k,l\in[n]$.
Then the {\em partition function} $p_R:\CC\to\oC$ is given by
\dyy{4se14g}{
p_R(C):=\sum_{\varphi:A\to[n]}\prod_{uv\in E}
R_{\varphi(u_{\text{in}}),\varphi(v_{\text{in}})}^{\varphi(u_{\text{out}}),\varphi(v_{\text{out}})}
}
for any multiloop chord diagram $C$,
where $A$ and $E$ denote the sets of directed and undirected edges, respectively, of $C$,
and where $v_{\text in}$ and $v_{\text{out}}$ denote the ingoing and the outgoing directed edge, respectively,
at a vertex $v$.
This implies $p_R(\loop)=n$.
Note that \rf{4se14g} is independent of the basis of $\oC^n$ chosen.

We will also write $p(C)(R)$ for $p_R(C)$.
Then $p(C):S^2(\gl(n))\to\oC$ is $\GL(n)$-invariant.
(Throughout, $\GL(n)$ acts on $\gl(n)$ by $h\cdot M:=hMh^{-1}$ for $h\in\GL(n)$ and $M\in\gl(n)$.)
By the First Fundamental Theorem (FFT) of invariant theory
(cf.\ [7] Corollary 5.3.2), each $\GL(n)$-invariant regular function
$S^2(\gl(n))\to\oC$ is a linear combination of functions $p(C)$ with $C$ a multiloop chord diagram.
(Here multiloop is essential.)

It will be convenient to notice at this point the following alternative
description of the partition function $p_R$.
Let $b_1,\ldots,b_k\in\gl(n)$ be as in \rf{4se14f}, with $X:=\gl(n)$.
Let $C$ be a multiloop chord diagram.
Consider a function $\psi:E\to[k]$.
`Assign' matrix $b_{\psi(uv)}$ to each of the ends $u$ and $v$ of any undirected edge $uv$.
Each of the Wilson loops in $C$ now has matrices assigned to its vertices, and on
each Wilson loop, we can take the trace of the product of these matrices (in order).
Taking the product of these traces over all Wilson loops,
and next summing up these products over all $\psi:E\to[k]$, gives $p_R(C)$.
(In the idiom of Szegedy [14], we here color the {\em undirected} edges, with $k$ colors,
while in \rf{4se14g} we color the {\em directed} edges, with $n$ colors.)

\tussenkop{Tangles}
We need an extension of the concept of multiloop chord diagram.
Define a {\em multiloop chord tangle}, or {\em tangle} for short,
as a graph with directed and undirected edges, such that
each vertex $v$ either satisfies \rf{22me14a} or $v$ is incident
with precisely one directed edge and with no undirected edge.
Of the latter type of vertex, there are two kinds:
vertices, called {\em roots}, with one outgoing edge, and
vertices, called {\em sinks}, with one ingoing edge.
The numbers of roots and of sinks are necessarily equal.
Again, a tangle may have components that are just
the vertexless directed loop \loop.

A {\em $k$-labeled multiloop chord tangle}, or just {\em $k$-tangle}, is a
tangle with precisely $k$ roots, equipped with labels $1,\ldots,k$, and
$k$ sinks, also equipped with labels $1,\ldots,k$.
Denote the collection of $k$-tangles by $\TT_k$.
So $\TT_0=\CC$.

For $S,T\in\TT_k$, let $S\cdot T$ be the multiloop chord
diagram arising from the disjoint union of $S$ and $T$ by,
for each $i=1,\ldots,k$, identifying the $i$-labeled sink in $S$
with the $i$-labeled root in $T$, and identifying the $i$-labeled root
in $S$ with the $i$-labeled sink in $T$;
after each identification, we ignore identified points
as vertex, joining its two incident directed edges into one directed edge;
that is,
\raisebox{-.3\height}{\scalebox{0.1}{\includegraphics{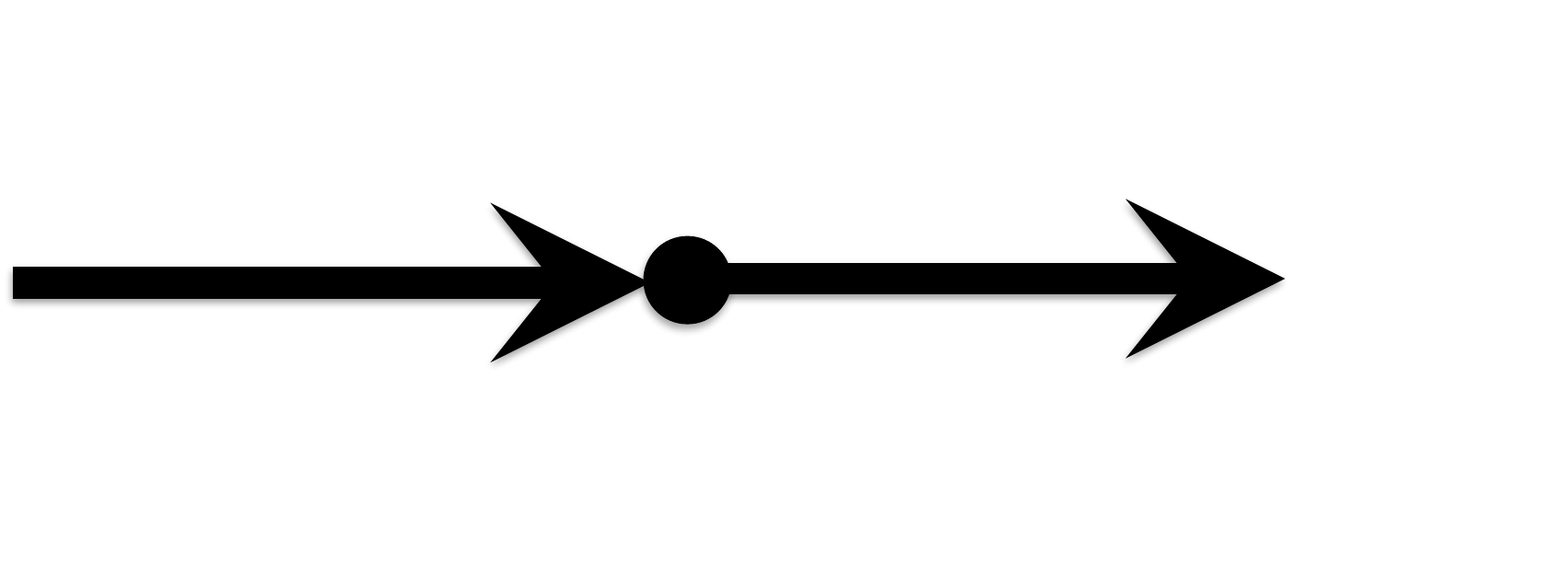}}}
becomes~~~
\raisebox{-.3\height}{\scalebox{0.1}{\includegraphics{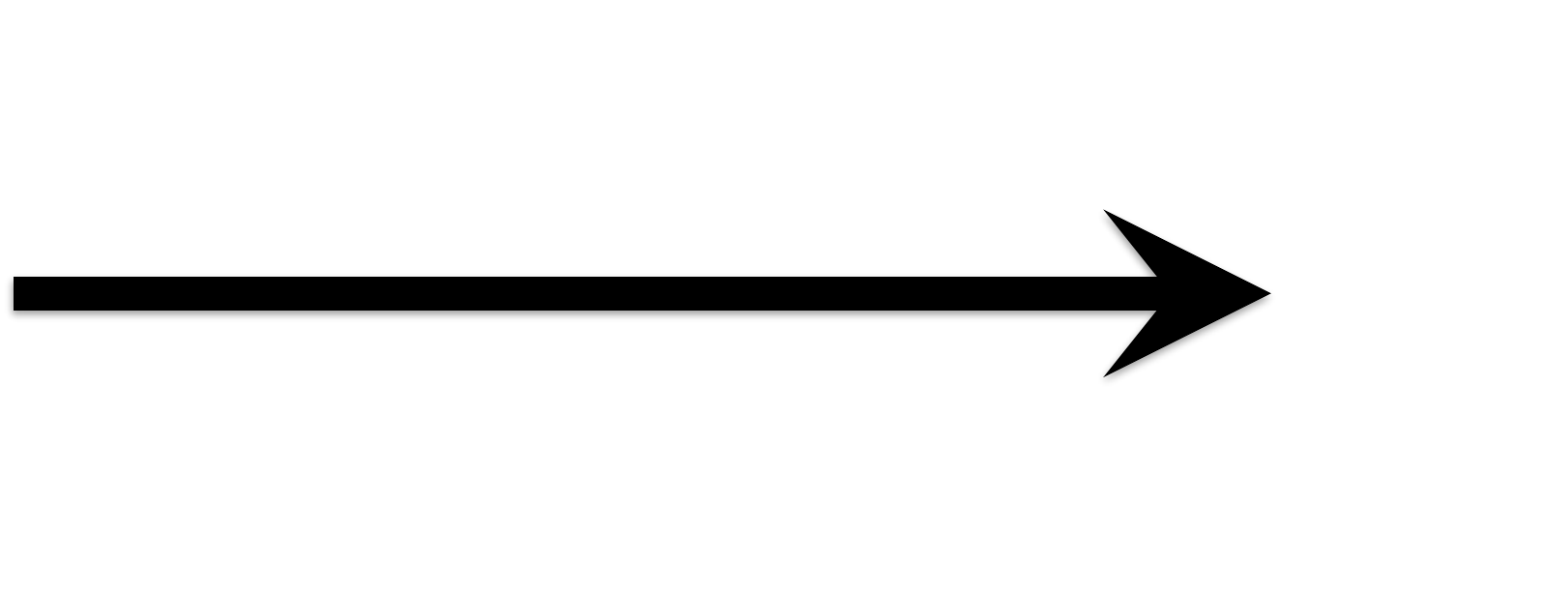}}}.
Note that this operation may introduce vertexless loops.
We extend this operation bilinearly to $\oC\TT_k$.
If $C,D\in\CC=\TT_0$, then $C\cdot D$ is equal to the disjoint union of $C$ and $D$.

Weight systems are determined by the 4T `quantum' 3-tangle $\tau_4$,
which is the element of $\oC\TT_3$ emerging from the 4T relations:
\dyy{3jn14b}{
\tau_4:=
\raisebox{-.3\height}{\scalebox{0.078}{\includegraphics{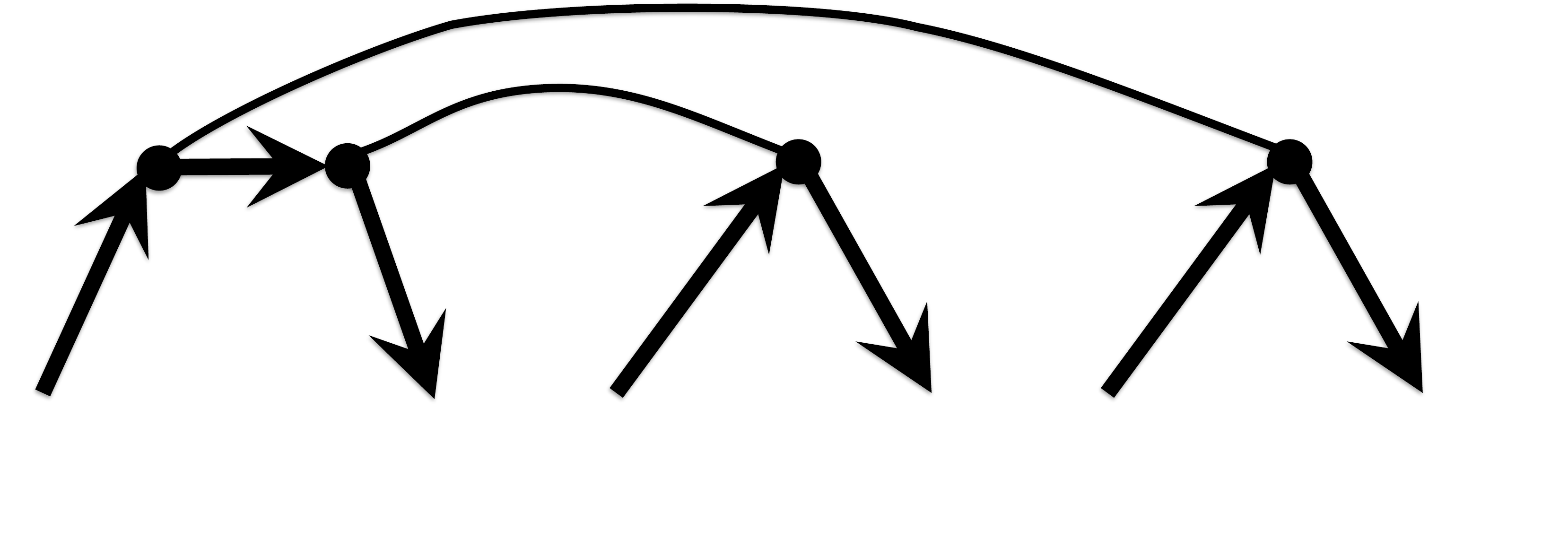}}}
-
\raisebox{-.3\height}{\scalebox{0.078}{\includegraphics{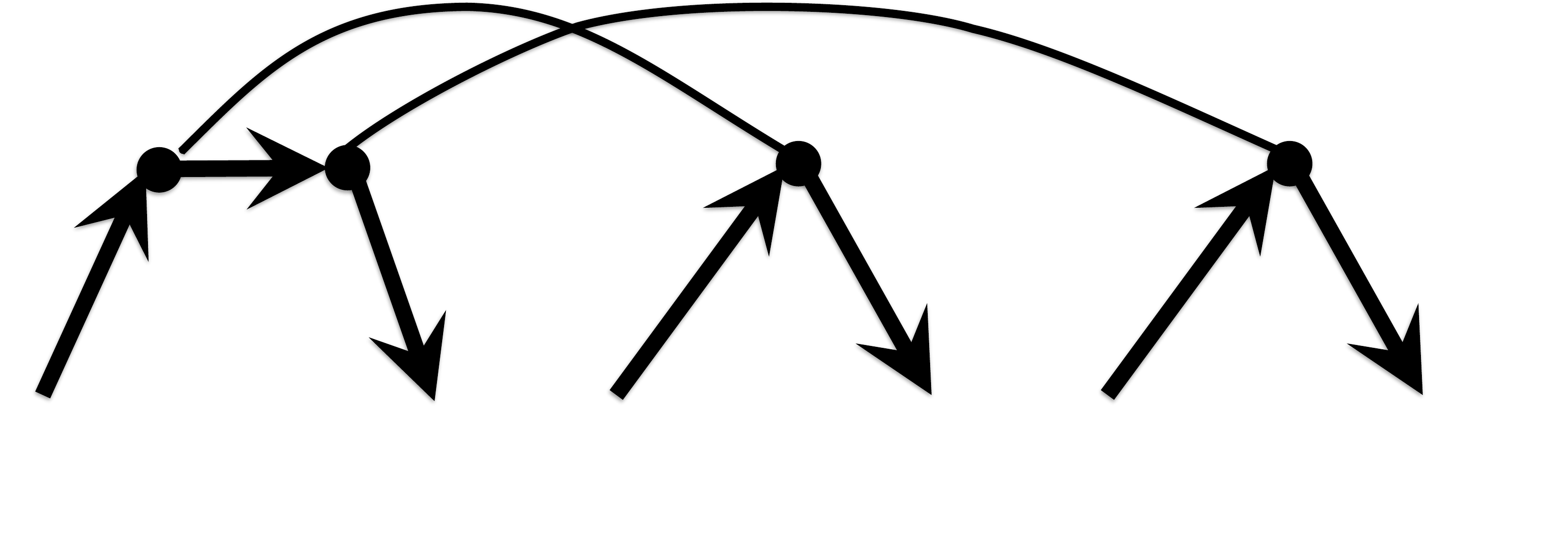}}}
-
\raisebox{-.3\height}{\scalebox{0.078}{\includegraphics{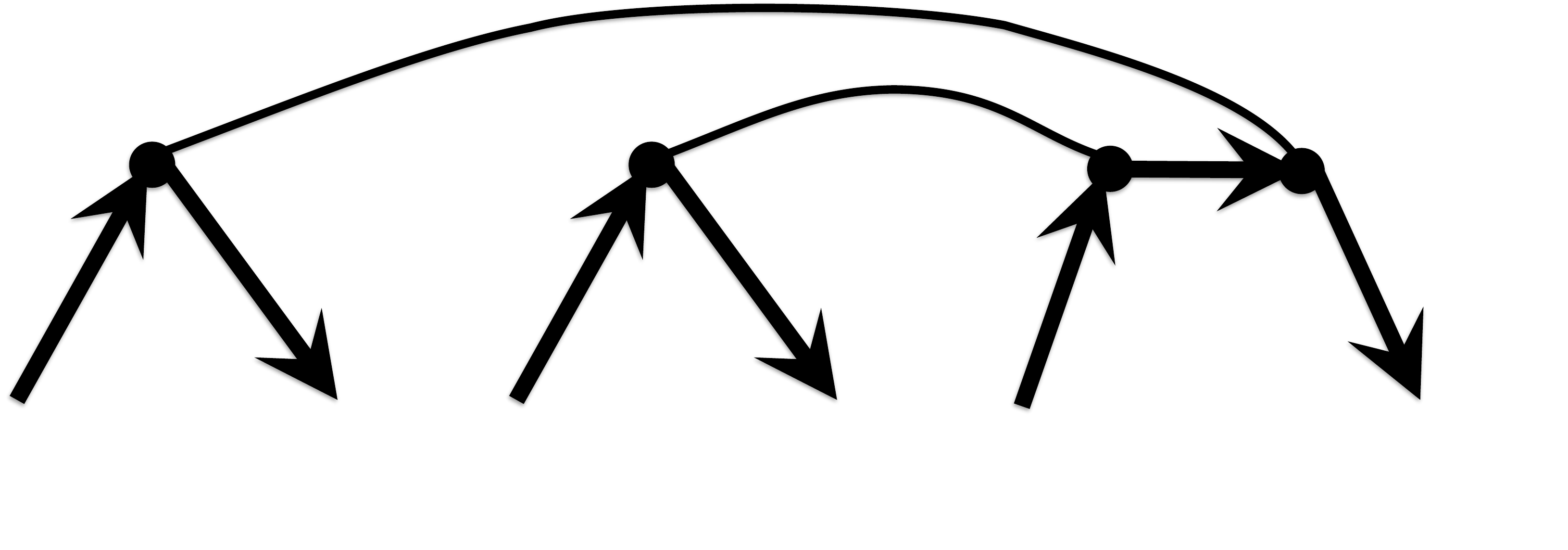}}}
+
\raisebox{-.3\height}{\scalebox{0.078}{\includegraphics{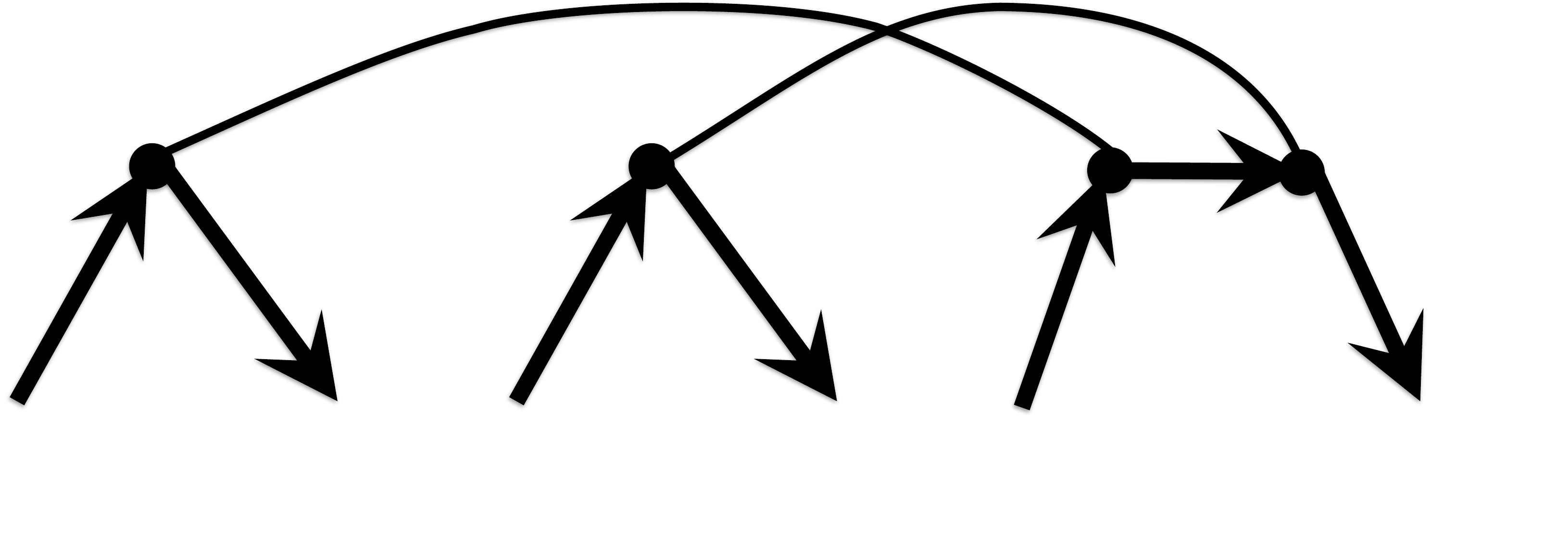}}}.
}
(We have omitted labels, as they are obvious (one may take labels $1,1,2,2,3,3$ from left to right in each
tangle in \rf{3jn14b}).)
Thus a function $f$ on $\CC$ is a weight system if and only if
$f(\tau_4\cdot T)=0$ for each 3-tangle $T$.

\tussenkop{The partition function on tangles}
We extend the function $p_R$ on multiloop chord diagrams to a function
$\widehat p_R$ on tangles.
For each $R\in S^2(\gl(n))$ and $k\in\oZ_+$,
the {\em partition function} $\widehat p_R:\TT_k\to \gl(n)^{\otimes k}$
is defined as, for $C\in\TT_k$:
\dyyz{
\widehat p_R(C):=\sum_{\varphi:A\to[n]}\prod_{uv\in E}
R_{\varphi(u_{\text{in}}),\varphi(v_{\text{in}})}^{\varphi(u_{\text{out}}),\varphi(v_{\text{out}})}
\bigotimes_{j=1}^kE_{\varphi(a_j)}^{\varphi(a_j^*)}.
}
Here we use the same notation as for \rf{4se14g}.
Moreover, $a_1,\ldots,a_k$ are the directed edges leaving the roots labeled $1,\ldots,k$,
respectively, and
$a_1^*,\ldots,a_k^*$ are the directed edges entering the sinks labeled $1,\ldots,k$, respectively.
For $h,i\in[n]$, $E_h^i$ is the matrix in $\gl(n)$ with 1 in position $(h,i)$ and 0 elsewhere.
Note that \rf{15me14a} is independent of the basis of $\oC^n$ chosen.

Again, set $\widehat p(C)(R):=\widehat p_R(C)$.
Then $\widehat p(C):S^2(\gl(n))\to\gl(n)^{\otimes k}$
is a $\GL(n)$-equivariant regular function,
and each such function is a linear combination of functions $\widehat p(C)$
(by the FFT for invariant theory).

Note that $p_R$ is the restriction of $\widehat p_R$ to $\CC$, and that
\dyy{14se14c}{
p_R(S\cdot T)=\tr(\widehat p_R(S)\widehat p_R(T))
}
for all $k$-tangles $S$ and $T$ (under the natural identification $\gl(n)^{\otimes k}=\End((\oC^n)^{\otimes k})$).

\tussenkop{Weight systems and Lie algebras}
A Lie algebra $\frak{g}$ is called {\em metrized} if it is equipped with a nondegenerate bilinear form $\langle.,.\rangle$
that is {\em ad-invariant}, i.e.\ satisfies $\langle[x,y],z\rangle=\langle x,[y,z]\rangle$
for all $x,y,z\in\frak{g}$.

For any $R\in S^2(\gl(n))$, choose linearly independent $b_1,\ldots,b_k\in\gl(n)$ such that
$R=\sum_{i=1}^kb_i\otimes b_i$ (as in \rf{4se14f}, taking $X:=\gl(n)$).
Then the following fundamental insight was given by Bar-Natan [1,\linebreak[0]2]:
\dy{13jn14b}{
$\widehat p_R(\tau_4)=0$
if and only if $b_1,\ldots,b_k$ form an orthonormal basis of a metrized Lie algebra $\frak{g}\subseteq\gl(n)$.
}
In fact, if $\frak{g}$ is a metrized Lie algebra and $\rho:\frak{g}\to\gl(n)$ is a representation,
then $R:=\sum_{i=1}^k\rho(b_i)\otimes\rho(b_i)$ satisfies $\widehat p_R(\tau_4)=0$
(where again $b_1,\ldots,b_k$ is any orthonormal basis of $\frak{g}$).
This implies that
\dyyz{
\varphi_{\frak{g}}^{\rho}:=p_R
}
is a weight system.

\sectz{Theorem and proof}

Define, for any $f:\CC\to\oC$ and $k\in\oZ_+$, the $\TT_k\times\TT_k$ matrix $M_{f,k}$ by
\dyyz{
(M_{f,k})_{S,T}:=f(S\cdot T),
}
for $S,T\in\TT_k$.

\thmnn{
Let $f:\CC\to\oC$ be a weight system.
Then the following are equivalent:
\di{15me14a}{
\nr{i} $f=\varphi_{\frak{g}}^{\rho}$ for some completely reducible faithful representation $\rho$ of some metrized Lie algebra $\frak{g}$;
\nrs{ii} $f=\varphi_{\frak{g}}^{\rho}$ for some representation $\rho$ of some metrized Lie algebra $\frak{g}$;
\nrs{iii} $f$ is the partition function $p_R$ of some $n\in\oZ_+$ and $R\in S^2(\gl(n))$;
\nrs{iv} $f(\loop)\in\oR$ and $\rank(M_{f,k})\leq f(\loop)^{2k}$ for each $k$.
}
}

\pf
(i){$\boldmath\Longrightarrow$}(ii) is trivial,
and
(ii){$\boldmath\Longrightarrow$}(iii) is easy by taking $R=\sum_{i=1}^k\rho(b_i)\otimes\rho(b_i)$ for
some orthonormal basis $b_1,\ldots,b_k$ of $\frak{g}$.

As to (iii){$\boldmath\Longrightarrow$}(iv):
$f(\loop)=n\in\oR$ is direct, while $\rank(M_{f,k})\leq n^{2k}$ follows from \rf{14se14c},
since $\widehat p_R(S)$ and $\widehat p_R(T)$ belong to $\gl(n)^{\otimes k}$, which is $n^{2k}$-dimensional.

It remains to show (iv)$\Longrightarrow$(i).
For $k\in\oZ_+$ and $S,T\in\TT_k$, define (next to the `inner product' $S\cdot T$)
the product $ST$ as the $k$-tangle obtained
from the disjoint union of $S$ and $T$ by identifying sink labeled $i$ of $S$ with root
labeled $i$ of $T$, and ignoring this vertex as vertex
(i.e.,
\raisebox{-.3\height}{\scalebox{0.1}{\includegraphics{becomes1.pdf}}}
becomes
\raisebox{-.3\height}{\scalebox{0.1}{\includegraphics{becomes2.pdf}}}),
for $i=1,\ldots,k$;
the roots of $S$ labeled $1,\ldots,k$ and sinks of $T$ labeled $1,\ldots,k$
make $ST$ to a $k$-tangle again.

Clearly, this product is associative, and satisfies $(ST)\cdot U=S\cdot(TU)$ for all $k$-tangles
$S,T,U$.
Moreover, there is a unit, denoted by ${\mathbf 1}_k$, consisting of $k$
disjoint directed edges $e_1,\ldots,e_k$, where both ends of $e_i$ are labeled $i$
($i=1,\ldots,k$).

Extend the product $ST$ bilinearly to $\oC\TT_k$,
making $\oC\TT_k$ to a $\oC$-algebra.
Let $\II_k$ be the null space of the matrix $M_{f,k}$, that is,
the space of $\tau\in\oC\TT_k$ with $f(\tau\cdot T)=0$
for each $k$-tangle $T$.
Then $\II_k$ is an ideal in the algebra $\oC\TT_k$, and the quotient
\dyyz{
\AAA_k:=\oC\TT_k/\II_k
}
is an algebra of dimension $\rank(M_{f,k})$.
We will indicate the elements of $\AAA_k$ just by their representatives
in $\oC\TT_k$.
Define the `trace-like' function $\vartheta:\AAA_k\to\oC$ by
\dyyz{
\vartheta(x):=f(x\cdot {\mathbf 1}_k)
}
for $x\in\AAA_k$.
Then $\vartheta(xy)=\vartheta(yx)$ for all $x,y\in\AAA_k$ and $\vartheta({\mathbf 1}_k)=f(\loop)^k=n^k$.

We first show that $\AAA_k$ is semisimple.
To this end, let for $k,m\in\oZ_+$ and $\pi\in S_m$,
$P_{k,\pi}$ be the $km$-tangle consisting of $km$ disjoint
edges $e_{i,j}$ for $i=1,\ldots,m$ and $j=1,\ldots,k$, where the head (sink) of $e_{i,j}$ is labeled $i+(j-1)m$
and its tail (root) is labeled $\pi(i)+(j-1)m$.

We also need a product $S\sqcup T$ of a $k$-tangle $S$ and an $l$-tangle $T$: it is the
$k+l$-tangle obtained from the disjoint union of $S$ and $T$ by adding $k$ to all labels in $T$.
This product can be extended bilinearly to $\oC\TT_k\times\oC\TT_l\to\oC\TT_{k+l}$.
The product is associative, so that for any $x\in\oC\TT_k$, the $m$-th power $x^{\sqcup m}$ is
well-defined.

Then for any $x\in\oC\TT_k$ and $\rho,\sigma\in S_m$ one has
\dyy{5se11a}{
f(x^{\sqcup m}P_{k,\rho}\cdot P_{k,\sigma})=
f(x^{\sqcup m}\cdot P_{k,\rho}P_{k,\sigma})=
f(x^{\sqcup m}\cdot P_{k,\rho\sigma})=
\prod_c\vartheta(x^{|c|}),
}
where $c$ ranges over the orbits of permutation $\rho\sigma$.
We are going to use that, for each $x\in\oC\TT_k$, the $S_m\times S_m$ matrix
$(f(x^{\sqcup m}P_{k,\rho}\cdot P_{k,\sigma}))_{\rho,\sigma\in S_m}$ has rank at most $\rank(M_{f,km})$
(since $x^{\sqcup m}P_{k,\rho}$ belongs to $\oC\TT_{km}$, for each $\rho$).

\cl{5se11b}{
For each $k$,
if $x$ is a nilpotent element of $\AAA_k$, then $\vartheta(x)=0$.
}

\pfcl
Suppose $\vartheta(x)\neq 0$ and $x$ is nilpotent.
Then there is a largest $t$ with $\vartheta(x^t)\neq 0$.
Let $y:=x^t$.
So $\vartheta(y)\neq 0$ and $\vartheta(y^s)=0$ for each $s\geq 2$.
By scaling we can assume that $\vartheta(y)=1$.

Choose $m$ with $m!>n^{2km}$.
By \rf{5se11a} we have, for any $\rho,\sigma\in S_m$,
\dyyz{
f(y^{\sqcup m}P_{k,\rho}\cdot P_{k,\sigma})=\delta_{\rho,\sigma^{-1}},
}
since $\vartheta(x^{|c|})=0$ if $|c|>1$, implying that the product in \rf{5se11a} is 0 if $\rho\sigma\neq\id$.

So $\rank(M_{f,km})\geq m!$, contradicting the fact that
$\rank(M_{f,km})\leq n^{2km}<m!$.
\obx

The following is a direct consequence of Claim \ref{5se11b}:

\clz{
$\AAA_k$ is semisimple, for each $k$.
}

\pfcl
As $\AAA_k$ is finite-dimensional, it suffices to
show that for each nonzero element $x$ of $\AAA_k$ there exists $y$ with $xy$ not nilpotent.
As $x\not\in\II_k$, we know that $f(x\cdot y)\neq 0$ for some $y\in \AAA_k$.
So $\vartheta(xy)\neq 0$, and hence, by Claim \ref{5se11b}, $xy$ is not nilpotent.
\obx

\cl{26jl11j}{
For each $k$, if $x$ is a nonzero idempotent in $\AAA_k$, then $\vartheta(x)$ is
a positive integer.
}

\pfcl
Let $x$ be any idempotent.
Then for each $m\in\oZ_+$ and $\rho,\sigma\in S_m$, by \rf{5se11a}:
\dyyz{
f(x^{\sqcup m}P_{k,\rho}\cdot P_{k,\sigma})=\vartheta(x)^{o(\rho\sigma)},
}
where $o(\pi)$ denotes the number of orbits of any $\pi\in S_m$.
So for each $m$:
\dyyz{
\rank((\vartheta(x)^{o(\rho\sigma)})_{\rho,\sigma\in S_m})\leq\rank(M_{f,km})\leq f(\loop)^{2km}.
}
This implies (cf.\ [13]) that
$\vartheta(x)\in\oZ$ and $\vartheta(x)\leq f(\loop)^k$.
As ${\mathbf 1}_k-x$ also is an idempotent in $\AAA_k$
and as $\vartheta({\mathbf 1}_k)=f(\loop)^k$,
we have $f(\loop)^k\geq\vartheta({\mathbf 1}_k-x)=f(\loop)^k-\vartheta(x)$.
So $\vartheta(x)\geq 0$.

Suppose finally that $x$ is nonzero while $\vartheta(x)=0$.
As $\vartheta(y)\geq 0$ for each idempotent $y$,
we may assume that $x$ is a minimal nonzero idempotent.
Let $J$ be the two-sided ideal generated by $x$.
As $\AAA_k$ is semisimple and $x$ is a minimal nonzero idempotent,
$J\cong\oC^{m\times m}$ for some $m$, yielding a trace function on $J$.
As $\vartheta$ is linear, there exists an $a\in J$ such that
$\vartheta(z)=\tr(za)$ for each $z\in J$.
As $\vartheta(yz)=\vartheta(zy)$ for all $y,z\in J$, we have $\tr(zay)=\tr(zya)$ for all $y,z\in J$.
So $ay=ya$ for all $y\in J$, hence $a$ is equal to
a scalar multiple of the $m\times m$ identity matrix in $J$.

As $x\neq 0$, $f(x\cdot y)\neq 0$ for some $y\in\AAA_k$, so $\vartheta(xy)\neq 0$.
Hence $a\neq 0$, and so $\vartheta(x)\neq 0$ (as $x$ is a nonzero idempotent), contradicting our assumption.
\obx

As ${\mathbf 1}_1$ is an idempotent in $\AAA_1$, Claim \ref{26jl11j} implies that
$f(\loop)=\vartheta({\mathbf 1}_1)$ is a nonnegative integer, say $n$.
Define an element $\Delta\in\oC\TT_{n+1}$ as follows.
For $\pi\in S_{n+1}$ let $T_{\pi}$ be the $(n+1)$-tangle
consisting of $n+1$ disjoint directed edges $e_1,\ldots,e_{n+1}$, where
the head of $e_i$ is labeled $i$ and its tail is labeled $\pi(i)$,
for $i=1,\ldots,n+1$.
Then
\dyyz{
\Delta:=\sum_{\pi\in S_{n+1}}\sgn(\pi)T_{\pi}.
}
Then $(n+1)!^{-1}\Delta$ is an idempotent in $\oC\TT_{n+1}$,
and
\dyyz{
\vartheta(\Delta)
=
\sum_{\pi\in S_{n+1}}\sgn(\pi)n^{o(\pi)}
=
\sum_{\pi\in S_{n+1}}\sgn(\pi)\sum_{\varphi:[n+1]\to[n]\atop\varphi\circ\pi=\varphi}1
=
\sum_{\varphi:[n+1]\to[n]}
\sum_{\pi\in S_{n+1}\atop\varphi\circ\pi=\varphi}\sgn(\pi)
=
0,
}
since $f(\loop)=n$ and since no $\varphi:[n+1]\to[n]$ is injective.
So by Claim \ref{26jl11j}, $\Delta=0$ in $\AAA_{n+1}$.
That is, $\Delta\in\II_{n+1}$.
So, by definition of $\II_{n+1}$,
\dyy{14se14b}{
\Delta\cdot\oC\TT_{n+1}\subseteq\Ker(f).
}
To conclude the proof of (iv)$\Longrightarrow$(iii), we follow a line of arguments similar to that in
[5].
Recall that $p:\oC\CC\to\OO(S^2(\gl(n)))$ is defined by $p(C)(X):=p_X(C)$ for all $C\in\CC$ and $X\in S^2(\gl(n))$.

\cl{12me14d}{
$\Ker~p\subseteq\Delta\cdot\oC\TT_{n+1}$.
}

\pfcl
Let $\gamma\in\oC\CC$ with $p(\gamma)=0$.
We prove that $\gamma\in\Delta\cdot\oC\TT_{n+1}$.
As each homogeneous component of $p(\gamma)$ is 0,
we can assume that $\gamma$ is a linear combination of multiloop chord diagrams that all have the same number $m$
of chords.

Let $H$ be the group of permutations of $[2m]$ that maintain the collection
$\{\{2i-1,2i\}\mid i\in[m]\}$.
Then $H$ naturally acts on $S_{2m}$ by $\rho\cdot\pi:=\rho\pi\rho^{-1}$ for $\rho\in H$
and $\pi\in S_{2m}$.

For any $\pi\in S_{2m}$, let $C_{\pi}$ be the multiloop chord diagram with with
vertex set $[2m]$,  chords $\{2i-1,2i\}$
(for $i=1,\ldots,m$) and directed edges $(i,\pi(i))$ (for $i=1,\ldots,2m$).
Each multiloop chord diagram with $m$ chords is isomorphic to $C_{\pi}$ for some $\pi\in S_{2m}$.
Therefore, we can write $\gamma=\sum_{\pi\in S_{2m}}\lambda(\pi)C_{\pi}$ with $\lambda:S_{2m}\to\oC$.
As $H$ leaves any $C_{\pi}$ invariant up to isomorphism, we can assume that $\lambda$ is $H$-invariant.

Define linear functions $F_{\pi}$ (for $\pi\in S_{2m}$) and $F$ on $\gl(n)^{\otimes 2m}$ by
\dyy{7se14c}{
F_{\pi}((a_1\otimes b_1)\otimes\cdots\otimes(a_{2m}\otimes b_{2m})):=\prod_{i=1}^{2m}b_i(a_{\pi(i)})
\text{ and }
F:=\sum_{\pi\in S_{2m}}\lambda(\pi)F_{\pi},
}
for $a_1,\ldots,a_{2m}\in\oC^n$ and $b_1,\ldots,b_{2m}\in(\oC^n)^*$.
Note that $F_{\pi}(R^{\otimes m})=p(C_{\pi})(R)$ for any $R\in S^2(\gl(n))$.
Hence $F(R^{\otimes m})=p(\gamma)(R)=0$.
We show that this implies that $F=0$.

Indeed, suppose $F(v_1\otimes\cdots\otimes v_{2m})\neq 0$ for some $v_1,\ldots,v_{2m}\in\gl(n)$.
For each $x\in\oC^m$, define
\dyyz{
R_x:=\sum_{i=1}^mx_i(v_{2i-1}\otimes v_{2i}+v_{2i}\otimes v_{2i-1})\in S^2(\gl(n)).
}
As $F$ is $H$-invariant (since $\lambda$ is $H$-invariant),
the coefficient of $x_1\cdots x_m$ in the polynomial $F(R_x^{\otimes m})$ is equal
to $|H|\cdot F(v_1\otimes\cdots\otimes v_{2m})\neq 0$.
So the polynomial is nonzero, hence $F(R_x^{\otimes m})\neq 0$ for some $x$, a contradiction.
Therefore, $F=0$.

Next, define the following polynomials $q_{\pi}$ (for $\pi\in S_{2m}$) and $q$ on $\oC^{2m\times 2m}$:
\dyy{7se14b}{
q_{\pi}(X):=\prod_{i=1}^{2m}X_{i,\pi(i)}
\text{ and }
q:=\sum_{\pi\in S_{2m}}\lambda(\pi)q_{\pi},
}
for $X=(X_{i,j})\in\oC^{2m\times 2m}$.
Then
\dy{6se14a}{
$q(X)=0$ if $\rank(X)\leq n$.
}
Indeed, if $\rank(X)\leq n$, then $X=(b_i(a_j))_{i,j=1}^{2m}$
for some $a_1,\ldots,a_{2m}\in\oC^n$ and $b_1,\ldots,b_{2m}\in(\oC^n)^*$.
By \rf{7se14c} and \rf{7se14b},
$q((b_i(a_j))_{i,j=1}^{2m})=F((a_1\otimes b_1)\otimes\cdots\otimes (a_{2m}\otimes b_{2m}))=0$,
proving \rf{6se14a}.

By the Second Fundamental Theorem (SFT) of invariant theory for $\GL(n)$
(cf.\ [7] Theorem 12.2.12),
\rf{6se14a} implies that $q$ belongs to the ideal in $\OO(\oC^{2m\times 2m})$ generated by the
${(n+1)}\times {(n+1)}$ minors of $\oC^{2m\times 2m}$.
Since each monomial in $q(X)$ contains precisely one variable from each row of $X$ and precisely one variable from each column of $X$, this implies $\gamma\in\Delta\cdot\oC\TT_{n+1}$.
\openbx

Claim \ref{12me14d} and \rf{14se14b} imply $\Ker(p)\subseteq\Ker(f)$, and so
there exists a linear function $\varphi:p(\oC\CC)\to\oC$ such that $\varphi\circ p=f$.
Then $\varphi$ is an algebra homomorphism, since for $C,D\in\CC$ one has
$\varphi(p(C)p(D))=\varphi(p(CD))=f(CD)=f(C)f(D)=\varphi(p(C))\varphi(p(D))$.

We now apply some more invariant theory.
As before, $\GL(n)$ acts on $\gl(n)$ by $h\cdot M:=hMh^{-1}$ for $h\in\GL(n)$ and $M\in\gl(n)$.
This action transfers naturally to $S^2(\gl(n))$.
By the FFT of invariant theory, we have
\dyy{15se14b}{
\OO(S^2(\gl(n)))^{\GL(n)}=p(\oC\CC).
}
So $\varphi$ is an algebra homomorphism $\OO(S^2(\gl(n)))^{\GL(n)}\to\oC$.
Hence the affine $\GL(n)$-variety
\dy{15se14c}{
$\VV:=\{R\in S^2(\gl(n))\mid q(R)=\varphi(q)$ for each $q\in\OO(S^2(\gl(n)))^{\GL(n)}\}$
}
is nonempty (as $\GL(n)$ is reductive).
By \rf{15se14b} and by substituting $q=p(C)$ in \rf{15se14c},
\dyyz{
\VV:=\{R\in S^2(\gl(n))\mid p_R=f\}.
}
Hence as $\VV\neq\emptyset$ we have \rf{15me14a}(iii).

To show \rf{15me14a}(ii), by \rf{13jn14b} it suffices to show that the function
$\widehat p(\tau_4):S^2(\gl(n))\to\gl(n)^{\otimes 3}$ has a zero on $\VV$.
Suppose $\widehat p(\tau_4)$ has no zero in $\VV$.
Then by the Nullstellensatz there exists a regular function $q:\VV\to\gl(n)^{\otimes 3}$
such that $\tr(\widehat p(\tau_4)(R)q(R))=1$ for each $R\in\VV$.
Applying the Reynolds operator, we can assume, as $\widehat p(\tau_4)$ is $\GL(n)$-equivariant, that
also $q$ is $\GL(n)$-equivariant.
Then by the FFT of invariant theory,
$q=\widehat p(\tau)$ for some $\tau\in\oC\TT_3$.
This gives $1= \tr(\widehat p(\tau_4)\widehat p(\tau)) = p(\tau_4\cdot\tau)$.
However, by \rf{15me14a}(iii), $p(\tau_4\cdot\tau)(R)=f(\tau_4\cdot\tau)=0$ for some $R$,
a contradiction.
This proves \rf{15me14a}(ii).

Finally, to show \rf{15me14a}(i), choose $R$ in the (unique) closed $\GL(n)$-orbit contained in $\VV$
(cf.\ [3], [11]).
Then $\widehat p_R(\tau_4)=0$, since by (ii), $\VV$ contains some $R'$ with
$\widehat p_{R'}(\tau_4)=0$, and since $R$ belongs to the closed orbit.
As $R$ belongs to $S^2(\gl(n))$, we can write
\dyy{16se14a}{
R=\sum_{i=1}^kb_i\otimes b_i
}
for some linearly independent $b_1,\ldots,b_k\in\gl(n)$,.
By \rf{13jn14b}, $b_1,\ldots,b_k$ form an orthonormal basis of a
metrized Lie subalgebra $\frak{g}$ of $\gl(n)$.

We prove that the identity $\id:{\frak{g}}\to\gl(n)$ is a completely reducible representation of $\frak{g}$.
Choose a chain $0=I_0\subset I_1\subset I_2\subset\cdots\subset I_{k-1}\subset I_k=\oC^n$ of
$\frak{g}$-submodules of $\oC^n$, with $k$ maximal.
For each $j=1,\ldots,k$, choose a subspace $X_j$ such that $I_j=I_{j-1}\oplus X_j$.
For each real $\lambda>0$, define $\Delta_{\lambda}\in\GL(n)$ by:
$\Delta_{\lambda}(x)=\lambda^jx$ if $x\in X_j$.

Then for each $M\in\frak{g}$, $M':=\lim_{\lambda\to \infty}\Delta_{\lambda}\cdot M$ exists.
Indeed, if $x\in X_j$, then $Mx\subseteq I_j$, and so
$\lim_{\lambda\to \infty}\Delta_{\lambda}M\Delta^{-1}_{\lambda}x$ is equal
to the projection of $Mx$ on $X_j$, with respect to the decomposition
$X_1\oplus\cdots\oplus X_k$ of $\oC^n$.
So $M'X_j\subseteq X_j$ for all $j$.

Hence, by \rf{16se14a}, also $R':=\lim_{\lambda\to \infty}\Delta_{\lambda}\cdot R$ exists, and is equal
to $\sum_{i=1}^kb_i'\otimes b'_i$.
As $\GL(n)\cdot R$ is closed, there exists $h\in\GL(n)$ with
$h^{-1}\cdot R=R'$, i.e., $R=h\cdot R'$.
Hence $\frak{g}$ is spanned by $h\cdot b'_1,\ldots,h\cdot b'_k$.
Therefore, $\frak{g}=\{h\cdot M'\mid M\in\frak{g}\}$.
Now $(h\cdot M')hX_j=hM'X_j\subseteq hX_j$ for each $M\in\frak{g}$ and $j$.
So $MhX_j\subseteq hX_j$ for each $M\in\frak{g}$ and $j$.
Therefore, each $hX_j$ is a $\frak{g}$-submodule.
By the maximality of $k$, each $hX_j$ is irreducible, proving \rf{15me14a}(i).
\bx

\section*{References}\label{REF}
{\small
\begin{itemize}{}{
\setlength{\labelwidth}{4mm}
\setlength{\parsep}{0mm}
\setlength{\itemsep}{1mm}
\setlength{\leftmargin}{5mm}
\setlength{\labelsep}{1mm}
}
\item[\mbox{\rm[1]}] D. Bar-Natan, 
Weights of Feynman diagrams and the Vassiliev knot invariants, 1991,
see
\url{http://www.math.toronto.edu/~drorbn/papers}

\item[\mbox{\rm[2]}] D. Bar-Natan, 
On the Vassiliev knot invariants,
{\em Topology} 34 (1995) 423--472.

\item[\mbox{\rm[3]}] M. Brion, 
Introduction to actions of algebraic groups,
{\em Les cours du C.I.R.M.} 1 (2010) 1--22.

\item[\mbox{\rm[4]}] S. Chmutov, S. Duzhin, J. Mostovoy, 
{\em Introduction to Vassiliev Knot Invariants},
Cambridge University Press, Cambridge, 2012.

\item[\mbox{\rm[5]}] J. Draisma, D. Gijswijt, L. Lov\'asz, G. Regts, A. Schrijver, 
Characterizing partition functions of the vertex model,
{\em Journal of Algebra} 350 (2012) 197--206.

\item[\mbox{\rm[6]}] M.H. Freedman, L. Lov\'asz, A. Schrijver, 
Reflection positivity, rank connectivity, and homomorphisms of graphs,
{\em Journal of the American Mathematical Society} 20 (2007) 37--51.

\item[\mbox{\rm[7]}] R. Goodman, N.R. Wallach, 
{\em Symmetry, Representations, and Invariants},
Springer, Dordrecht, 2009.

\item[\mbox{\rm[8]}] P. de la Harpe, V.F.R. Jones, 
Graph invariants related to statistical mechanical models:
examples and problems,
{\em Journal of Combinatorial Theory, Series B} 57 (1993) 207--227.

\item[\mbox{\rm[9]}] V. Kodiyalam, K.N. Raghavan, 
Picture invariants and the isomorphism problem for complex semisimple Lie algebras,
2004,
ArXiv \url{http://arxiv.org/math/0402215v1}

\item[\mbox{\rm[10]}] M. Kontsevich, 
Vassiliev's knot invariants. I.M. Gelfand seminar Part 2,
{\em Advances in Soviet Mathematics} 16, American Mathematical Society,
Providence, R.I., 1993, pp. 137--150.

\item[\mbox{\rm[11]}] H. Kraft, 
{\em Geometrische Methoden in der Invariantentheorie},
Vieweg, Braunschweig, 1984.

\item[\mbox{\rm[12]}] L. Lov\'asz, 
{\em Large Networks and Graph Limits},
American Mathematical Society, Providence, R.I., 2012.

\item[\mbox{\rm[13]}] A. Schrijver, 
Characterizing partition functions of the vertex model by rank growth,
preprint, 2012,
ArXiv \url{http://arxiv.org/abs/1211.3561}

\item[\mbox{\rm[14]}] B. Szegedy, 
Edge coloring models and reflection positivity,
{\em Journal of the American Mathematical Society}
20 (2007) 969--988.

\end{itemize}
}

\end{document}